\documentclass[a4,12pt,onecolumn,oneside,final]{article}
\usepackage{amssymb}
\usepackage{amscd}

\def\beq{\begin{equation}}
\def\eeq{\end{equation}}
\def\bea{\begin{eqnarray}}
\def\eea{\end{eqnarray}}
\def\nn{\nonumber}
\def\h#1{\hat{#1}}
\def\p#1{\partial_{#1}}
\def\A{{\cal A}}
\def\hot{\frac{h}{2}}

\begin{document}
\setlength{\baselineskip}{18pt}
\begin{center}
\thispagestyle{empty}
  \vspace*{3cm}
  {\LARGE \sf Noncommutative Geometry of Super-Jordanian 
   $OSp_h(2/1)$ Covariant Quantum Space
  }
  
  \vspace{1cm}
   N. Aizawa \\ \vspace{3mm}
   \textit{Department of Applied Mathematics, \\
           Osaka Women's University, \\
           Sakai, Osaka 590-0035, Japan} \\
   \vspace{1cm}
   R. Chakrabarti\footnote{Permanent address : Department of Theoretical Physics, 
           University of Madras, 
           Guindy Campus, Chennai, 600025, India} \vspace{3mm}
   \\ \vspace{3mm}
   \textit{Institute of Mathematical Sciences, Chennai 600113, India}
\end{center}
\vfill
\begin{abstract}
   Extending a recently proposed procedure of construction of 
various elements of differential geometry on noncommutative 
algebras, we obtain these structures on noncommutative superalgebras. 
As an example, a quantum superspace covariant under the action of 
super-Jordanian $OSp_h(2/1)$ is studied. It is shown that there 
exist a two-parameter family of torsionless connections, and the 
curvature computed from this family of connections is bilinear. It 
is also shown that the connections are not compatible with the metric. 
\end{abstract}
\newpage
%
%
%
\setcounter{equation}{0}
\section{Introduction}
\label{Intro}

  Noncommutative geometry is one of the active fields in recent 
theoretical physics and mathematics. The physical interests seem to be 
focused on noncommutative differential geometry, since it plays 
crucial roles in the context of string theory and quantum gravity 
\cite{Con,SZ}. There exist various approaches to noncommutative 
differential calculus. For instance, Connes' approach involves the 
Dirac operator \cite{Con}, Dubois-Violette's method is based on 
derivations \cite{DV1}, and studies  based on quantum groups 
also exist \cite{Wor,Manin,WZ}. Analogues of Riemannian connection, 
curvature and metric on noncommutative algebra $ \A $ were 
introduced in Ref. \cite{CFF} where the authors used only 
the left $ \A$-module structure of the differential forms. 
On the other hand, the $\A$-bimodule structure of an algebra of 
differential forms was used to define a linear connection for a  
particular differential calculus based on derivations \cite{DVM}. 
Mourad also made essential use \cite{Mou} of $\A$-bimodule 
structure to define  linear connection, torsion and curvature on a       
noncommutative algebra $\A$. In the procedure used in \cite{Mou}, 
a noncommutative generalization of the permutation operator on two 
copies of one-forms was introduced. The generalized permutation 
operator plays a role in defining the noncommutative differential 
geometry. This methodology was found to be useful in other 
approaches to the noncommutative differential calculi.  
Furthermore, it was extended to other studies concerning differential 
calculi on noncommutative algebras such as $ SL_q(2) $ covariant 
quantum plane \cite{DVMMM}, two-parameter quantum plane \cite{GMW}, 
Jordanian $h$-deformed quantum plane \cite{KSA,CMP}, matrix geometries 
\cite{MMM}, and so on. The curvatures corresponding to the connections 
\cite{DVM,Mou} were also studied \cite{DVMMMcurv} in this context. 

  In the present work, we follow the line of investigation developed in 
Refs. \cite{DVM,Mou,DVMMMcurv} to study the noncommutative differential
geometry associated with the quantized supergroups. It will be seen that 
the ideas used there are also appropriate for studying noncommutative 
differential geometry on quantum superspaces. Following Refs. 
\cite{DVM, Mou, DVMMMcurv} we can naturally define linear connections, 
torsions, curvatures and metrics on noncommutative quantum  
superspaces. As an example, we will study the quantum superspace 
covariant under the action of super-Jordanian deformation of $ OSp(2/1).$ 
The super-Jordanian $OSp_h(2/1)$ is introduced as a Hopf algebra dually 
related to the recently obtained triangular deformation of Lie superalgebra 
$ osp(2/1)$ \cite{ACS}; and it coincides with the deformed $OSp(2/1)$ 
supergroup studied by Juszczak and Sobczyck \cite{JS}. We obtain the most 
general form of linear connections on the quantum superspace covariant
under the action of $OSp_h(2/1)$. The curvatures and the metric will also be 
studied in the sequel.

  We first briefly review the construction of a linear connection in 
a commutative geometry \cite{Mou}. Let $ {\cal M} $ be a manifold, and 
$ C({\cal M}) $ be an algebra of functions on ${\cal M}$. The set of 
$k$-forms is denoted by $ \Omega^k.$ The covariant derivative $D$ is a 
linear map from $\Omega^1$ to $\Omega^1 \otimes_{C({\cal M})} \Omega^1$ 
obeying the Leibnitz rule
\beq
  D(f \xi) = df \otimes \xi + f D \xi, \qquad
  f \in C({\cal M}),\quad \xi \in \Omega^1.
  \label{Lcomm1}
\eeq
Since functions commute with forms, the Leibnitz rule may also be 
recast as 
\beq
  D(\xi f) = \sigma(\xi \otimes df) + (D\xi) f,  \label{Lcomm2}
\eeq
where $ \sigma $ is a permutation acting on 
$ \Omega^1 \otimes_{C({\cal M})} \Omega^1 $ 
\beq
  \sigma(\xi \otimes \eta) = \eta \otimes \xi, \qquad 
  (\xi,\,\eta) \in \Omega^1
  \label{sigmacomm}
\eeq
and the exterior derivative $d$ is nilpotent: $d^2 = 0$.
In the context of commutative geometry the two Leibnitz rules 
(\ref{Lcomm1}) and (\ref{Lcomm2}) are equivalent. However, in the 
noncommutative setting this is not the case, as functions and forms 
do not commute. The noncommutative covariant derivative is constructed
in such a way that it is required to satisfy \cite{Mou} the two 
Leibnitz rules. Reflecting the noncommutative nature of functions and 
forms, the operator $ \sigma $ can no longer be represented by a  
simple permutation element. It needs to be modified for noncommutative
quantum spaces. 

  Suppose the manifold ${\cal M} $ is parallelizable and let $ \omega^i $ 
be an arbitrary basis element of $\Omega^1$. The covariant derivative 
of a one-form is then uniquely determined by $ D \omega^i. $ The linear 
connection is defined by $ \Gamma^i = -D \omega^i. $ Namely, 
the covariant derivative defines the linear connection. Throughout this 
article, we use the terms `linear connection' and `covariant derivative' 
synonymously. Let $\pi$ be a projection of 
$ \Omega^1 \otimes_{C({\cal M})} \Omega^1 $ onto 
$ \Omega^2 $ such that $ \pi(\xi \otimes \eta) = \xi \wedge \eta. $ 
Then the map $ \Theta : \Omega^1 \rightarrow \Omega^2 $ defined by 
$ \Theta = d - \pi \circ D $ is a bimodule homomorphism, that is, 
it maintains $ \Theta(f \xi) = f\Theta(\xi) $ and $ \Theta(\xi f) 
= \Theta(\xi) f$, where $f \in C({\cal M})$. 
The torsion is defined by $\Theta(\omega^i).$ 
This construction of linear connections will be extended to noncommutative 
superspaces associated with quantum supergroups.

 This paper is organized as follows: In the next section, 
we extend the differential geometry on noncommutative algebras to 
noncommutative superalgebras. The super-Jordanian deformation of 
$OSp(2/1)$ is introduced in Section \ref{SJSOP21}. The quantum 
superspace which is covariant under the action of super-Jordanian 
$ OSp_h(2/1)$ is introduced, and the differential calculus on it in 
the sense of Wess-Zumino is constructed in Section \ref{DiffCal}. 
The linear connection on the quantum superspace is studied in Section 
\ref{Connection}, and it is observed that the most general torsionless 
connection is a member of a two-parameter family. In Section 
\ref{Curv-Met}, the curvature obtained from the linear connection 
is calculated and it is shown that the curvature is bilinear. The 
metric of the quantum superspace is also studied. We show that the 
covariant derivative is not compatible with the metric. 
Section \ref{CR} contain the concluding remarks. 

%
%
%
\setcounter{equation}{0}
\section{Noncommutative extension of superspace geometry}
\label{Extension}

  Let $ \A $ be a noncommutative algebra with ${\mathbb Z}_2$ grading. 
The grading is specified by parity of elements of $ \A. $ 
An even (odd) element $f \in \A $ has a parity $ \hat{f} = 0\,(1)$.
It is assumed that a differential calculus over $ \A $, describing, 
in particular, the one-forms and their commutation relations 
with the elements of $ \A $, has been constructed. 
Let  $ \Omega^k(\A) $ and $d$ denote the space of $k$-forms over $\A$ 
and the exterior derivative, respectively. 
The covariant derivative $ D $ is defined as a map 
 $ D : \Omega^1 \rightarrow \Omega^1 \otimes_{\A} \Omega^1 $ 
subject to the following Leibnitz rules
 \bea
   & & D(f \xi) = df \otimes \xi + (-1)^{\hat f} f D\xi, \label{D1} \\
   & & D(\xi f) = (-1)^{\hat \xi} \sigma(\xi \otimes df) + (D \xi) f, \qquad
       f \in \Omega^0, \quad \xi \in \Omega^1, \label{D2}
 \eea
 where 
 $ \sigma : \Omega^1 \otimes_{\A} \Omega^1 \rightarrow \Omega^1 
\otimes_{\A} \Omega^1 $ 
refers to a noncommutative generalization of the permutation map. 
The covariant derivative changes the parity of a $k$-form 
$ \xi $ by unity: $ \widehat{D\xi} = \hat{\xi} + 1, $ (mod 2). 
In the commutative case where 
$ \sigma (\xi \otimes \eta ) = (-1)^{\hat \xi \hat \eta} \eta 
\otimes \xi $, the two Leibnitz rules (\ref{D1}) and (\ref{D2}) are 
equivalent up to an overall sign
\beq
  D(\xi f) = (-1)^{\hat{\xi}\hat{f}} D(f\xi).  \label{Dclassical}
\eeq 
Using the definition of the covariant derivative, one can show that 
the map $ \sigma $ is $ \A$-bilinear
\beq
   \sigma(f \xi \otimes \eta ) = f \sigma(\xi \otimes \eta), \qquad 
   \sigma(\xi \otimes \eta f) = \sigma(\xi \otimes \eta) f, 
   \qquad f \in \Omega^0,    \quad 
   (\xi, \eta) \in \Omega^1.  \label{sig-lin}
\eeq
We now demonstrate the first relation in (\ref{sig-lin}). The second 
relation in (\ref{sig-lin}) follows similarly. For arbitrary
elements $ (f, g) \in \Omega^0 $ and $ \xi \in \Omega^1, $ we compute 
$ D(f \xi g) $ in two different 
ways. Regarding it as $ D(f \cdot \xi g), $ 
we apply (\ref{D1}) and obtain
\[
   D(f \xi g) = df \otimes \xi g + (-1)^{\hat f} f D(\xi g) 
   =  df \otimes \xi g + (-1)^{\hat f + \hat{\xi}} f \sigma(\xi \otimes dg) 
   + (-1)^{\hat f} f(D\xi)g. 
\]
Alternately, for the choice  $ D(f\xi g) = D(f\xi \cdot g)$ the 
Leibnitz rule (\ref{D2}) yields
\[
   D(f \xi g) 
   = (-1)^{\hat f + \hat{\xi}} \sigma(f\xi \otimes dg) + df \otimes \xi g + 
   (-1)^{\hat f} f(D\xi)g. 
\]
As the above two computations must give identical results, it follows 
$ \sigma(f\xi \otimes dg) = f\sigma(\xi \otimes dg). $ 
 
  The covariant derivative may be extended as a linear map from the 
$n$-fold tensored  space $ \otimes^n \Omega^1 $ to the $(n+1)$-fold 
tensored space $ \otimes^{(n+1)} \Omega^1. $ This is done recurrently 
while maintaining the following extension of the Leibnitz rule
\beq
   D (\omega \otimes \omega') = D \omega \otimes \omega' 
   + (-1)^{\hat{\omega}}\sigma_{12} ( \omega \otimes D \omega'),
   \label{Donn}
\eeq
where $ \omega \in \Omega^1, \ \omega' \in \otimes^{n-1} \Omega^1 $ and 
$ \sigma_{12} $ has a nontrivial structure in the first two sectors:
\beq
  \sigma_{12} = \sigma \otimes \underbrace{1 \otimes 1 \otimes \cdots 
  \otimes 1}_{n-1\   {\rm times}}.
  \label{sigma-s}
\eeq
Let $ \pi $ be a projection of $ \Omega^1 \otimes_{\A} \Omega^1 $ onto 
$ \Omega^2 $ defined by  the wedge product on the forms
\beq
    \pi(\xi \otimes \eta) = \xi \wedge \eta. \label{pidef}
\eeq
The noncommutativity of $ \A $, in general, demands 
$ \xi \wedge \eta \neq -\eta \wedge \xi. $  
Employing the projection operator $ \pi, $ we define the torsion 
$ \Theta $ of the covariant derivative $D$ 
as a map $ \Theta : \Omega^1 \rightarrow \Omega^2 $
\beq
   \Theta : \Omega^1 \rightarrow \Omega^2 , \qquad 
   \Theta = d - \pi \circ D.        \label{torsion}
\eeq 
The torsion is always left $\A$-linear, whereas the condition
\beq
     \pi \circ (\sigma - 1) = 0 \label{pisigma}
\eeq
is necessary for it to be right $\A$-linear.  
More explicitly, the torsion satisfies the relations
\beq
     \Theta(f \xi) = (-1)^{\hat f} f\Theta(\xi), \qquad
     \Theta(\xi f) = \Theta(\xi) f, \qquad 
     f \in \Omega^0,\ \xi \in \Omega^1. \label{torsion-lin}
\eeq
The condition (\ref{pisigma}) is necessary for the validity of the 
second relation in (\ref{torsion-lin}). 
Note that the relation (\ref{pisigma}) 
has a sign difference from the nongraded case \cite{DVMMM,MMM}. 
Since the proof is straightforward, we show only the second 
relation. The exterior derivative acts on the one-form $ \xi f$
as follows
\[
  d(\xi f) = (d\xi) f + (-1)^{\hat{\xi}} \xi \wedge df
  = (d\xi)f + (-1)^{\hat{\xi}}\, \pi(\xi \otimes df),
\] 
while the action of $ \pi \circ D $ on $ \xi f $ reads
\[
  \pi \circ D (\xi f) = \pi ( (-1)^{\hat{\xi}}\, 
  \sigma(\xi \otimes df) + (D\xi) f).
\]
Consequently, it follows
\[
  \Theta (\xi f) = \Theta(\xi) f - (-1)^{\hat{\xi}}\, 
  \pi \circ (\sigma-1) (\xi \otimes df).
\]
It is thus evident that the condition (\ref{pisigma}) need to be 
satisfied for the torsion $\Theta$ to be right $\A$-linear. 

   The curvature is defined by the following map \cite{DVMMMcurv}
\beq
    \pi_{12} D^2 \; : \; \Omega^1 \rightarrow \Omega^2 \otimes_{\A} \Omega^1,
    \label{curvdef}
\eeq
where $ \pi_{12} = \pi \otimes 1. $ The torsionless condition 
$ \Theta = 0 $ and the validity of the constraint (\ref{pisigma})  
requires the curvature to be  left $\A$-linear:
\beq
    \pi_{12}D^2(f\xi) = f\pi_{12}D^2(\xi), \qquad f \in \Omega^0, \ 
    \xi \in \Omega^1.
    \label{curvlinear}
\eeq
We demonstrate this below. Employing the Leibnitz rule (\ref{Donn}), 
we compute
\[
    D^2(f\xi) = Ddf \otimes \xi + (-1)^{\widehat{df}} 
    \sigma_{12}(df \otimes D\xi) 
    + (-1)^{\hat f} df \otimes D\xi + f D^2 \xi.
\]
The left hand side in (\ref{curvlinear}) now reads
\[
    \pi_{12} D^2(f\xi) = \pi\circ Ddf \otimes \xi - (-1)^{\hat f}\pi_{12} \circ 
    (\sigma_{12}-1)(df \otimes D\xi)
    + f \pi_{12} D^2 (\xi).
\]
The first and the second terms in the above expression vanish because of 
the torsionless condition and the constraint (\ref{pisigma}), respectively. 
Thus the curvature is left $ \A$-linear. In general, the curvature 
is not right $\A$-linear. It is, however, known that there exist some 
cases for nongraded $\A$ where the curvature is right $\A$-linear 
\cite{DVMMM}. We will find such an example for graded $\A$ in the 
following sections. 

  Now let us define a metric. A metric $g$ is a non-degenerate $\A$-bilinear 
map 
\beq
   g \; : \; \Omega^1 \otimes_{\A} \Omega^1 \rightarrow \A.
   \label{metricdef}
\eeq
The metric is said to be non-degenerate if the following conditions 
hold: $ g(\xi \otimes \eta) = 0$ for all $\eta \in \Omega^1$ implies 
$ \xi = 0, $ and, simultaneously, $ g(\xi \otimes \eta) = 0 $ for 
all $ \xi \in \Omega^1 $ implies $ \eta = 0. $ Symmetry of the metric 
is defined by using the extended permutation $ \sigma. $ A metric 
satisfying $ g \circ \sigma = g\, (g \circ \sigma = -g)$ is known to be 
symmetric (skew-symmetric) in nature. If the following diagram is
commutative, the covariant derivative $D$ is said 
to be compatible with the metric $g$, or, in short, 
$ D $ is said to be metric: 
\[
    \begin{CD}
    \Omega^1 \otimes_{\A} \Omega^1 @> D >> \Omega^1 \otimes_{\A} 
    \Omega^1 \otimes_{\A} \Omega^1 \\
    @V g VV   @V 1 \otimes g VV \\
    \A @> d >> \Omega^1
    \end{CD}
\]
More explicitly, the above compatibility condition reads
\beq
    d \circ g = (1 \otimes g) \circ D.   \label{Dtog}
\eeq

  In the following sections an example of the differential geometry 
described here will be presented. In this example, the algebra 
$ \A $ is taken to be a quantum superspace covariant under the 
action of a quantum supergroup $ OSp_h(2/1). $ The example 
will be constructed so as to keep the covariance of all relations. 

%
%
%
\setcounter{equation}{0}
\section{Super-Jordanian deformation of $ OSp(2/1) $}
\label{SJSOP21}

In this Section we introduce a quantum deformation of the supergroup 
$OSp(2/1)$. The conventions adopted here regarding the graded 
Yang-Baxter equation are same as in the Refs. \cite{KS,KR}. 
The quantum supergroup discussed here is the dual Hopf algebra 
to the super-Jordanian deformed $U_h(osp(2/1))$ algebra 
introduced recently. The study of super-Jordanian $ osp(2/1) $ 
algebra was initiated by Kulish \cite{Ku}. It was further developed 
by the works of the present authors \cite{ACS} and Borowiec 
$et. \ al.$ \cite{BLT}. In Ref.\cite{ACS}, the universal ${\cal R}$ matrix 
of the $U_h(osp(2/1)) $ algebra was obtained up to $ O(h^3) $ where $h$ 
is the deformation parameter. Its limiting classical value is 
described by $h \to 0$. The fundamental representation of the 
generators of the $ U_h(osp(2/1)) $ algebra is obtained by mapping the 
deformed algebra on its classical counterpart. Although the two 
deformation maps given in Ref. \cite{ACS} provide two 
distinct sets of matrices for the fundamental representation, 
the pertinent $R$ matrices computed for these two cases are identical. 
All the terms in the universal ${\cal R}$ matrix $O(h^3)$ and above  
vanish in the fundamental representation, and, therefore, the 
$R$ matrix in the said representation is determined by the terms up to 
$O(h^2). $ 
The $R$ matrix, thus obtained, is given by
\beq
  R =  
    \left(
      \begin{array}{ccc|ccc|ccc}
         1 & \cdot & -h & \cdot & h &  \cdot & h & \cdot & h^2/2 \\
         \cdot & 1 &  \cdot & \cdot & \cdot & -h & \cdot & \cdot &  \cdot \\
         \cdot & \cdot &  1 & \cdot & \cdot &  \cdot & \cdot & \cdot & -h \\ 
         \hline
         \cdot & \cdot & \cdot & 1 & \cdot &  \cdot & \cdot & h &  \cdot \\
         \cdot & \cdot & \cdot & \cdot & 1 &  \cdot & \cdot & \cdot & -h \\
         \cdot & \cdot & \cdot & \cdot & \cdot &  1 & \cdot & \cdot &  \cdot \\ 
         \hline
         \cdot & \cdot & \cdot & \cdot & \cdot & \cdot & 1 & \cdot &  h \\
         \cdot & \cdot & \cdot & \cdot & \cdot & \cdot & \cdot & 1 &  \cdot \\
         \cdot & \cdot & \cdot & \cdot & \cdot & \cdot & \cdot & \cdot &  1
      \end{array}
    \right)_,    \label{Rinfun}
\eeq
where the dot $(\cdot)$ is used instead of 0 for better readability. 
The $R$ matrix (\ref{Rinfun}) solves the graded Yang-Baxter equation. 
The inverse of this $R$ matrix is given by $ R^{-1} = R(-h) $; and 
identifying $h = -p$ in (\ref{Rinfun}) the $R$-matrix given in 
Ref. \cite{JS} is reproduced. In Ref. \cite{Ku}, a contraction technique 
is applied to the $R$ matrix in the fundamental representation of the 
$ U_q(osp(2/1)) $ algebra to obtain a triangular $\tilde R$ matrix , 
which  maintains the relation $\tilde R^{k\ell}_{ij} = R^{\ell k}_{ji}$ 
with the $R$ matrix given in (\ref{Rinfun}).  

  Now we explicitly write down the nonstandard deformed supergroup 
$OSp_h(2/1)$. Since the $R$ matrix (\ref{Rinfun}) is the 
inverse of the one used in Ref. \cite{JS}, the quantum supergroup 
$OSp_h(2/1)$ is identical to the one given in Ref. \cite{JS} where the 
deformed supergroup $OSp_h(2/1)$ is 
constructed by the FRT \cite{FRT} method. Let the inverse scattering 
matrix $T$ in the fundamental representation of the super-Jordanian 
deformed $OSp_h(2/1)$ is given by 
\beq
  T = (t^i_j) = \left(
   \begin{array}{ccc}
      a & \alpha & b \\
      \gamma & e & \beta \\
      c & \delta & d
   \end{array}
  \right)_,           \label{Osp}
\eeq
where $ \hat{i} = 0 (1)$ for $ i = \{1, 3\} (\{2\}) $, and 
$ \hat{t}^i_j = \hat{i}+\hat{j}. $ Thus the entries 
$ a, b, c, d, e $ are even elements and $ \alpha, \beta, \gamma, \delta $ 
are odd ones. The RTT relation and deformed orthosymplectic conditions
\footnote{The matrix $J$ corresponds to the matrix $C^{-1}$ in \cite{JS}.
An error contained in $C$ of \cite{JS} is corrected.}
\beq
  T^{st} J T = J, \qquad T J^{-1} T^{st} = J^{-1}, \qquad
  J = \left(
    \begin{array}{ccc}
      0 &  0 & 1 \\
      0 &  1 & 0 \\
     -1 &  0 & -h/2
    \end{array}
   \right)   \label{Osymh}
\eeq
determine the relations among the entries of $T.$  The 
supertransposition of $T$ is defined by $(T^{st})^i_j 
= (-1)^{\hat{i}(\hat{i}+\hat{j})} t^j_i$. 
It follows from this definition that $ (AB)^{st} = B^{st} A^{st}, 
\ ((A^{st})^{st})^i_j = (-1)^{\hat{i}+\hat{j}} A^i_j. $ 
Note that the matrix $J$ has the property
\beq
  (-1)^{\hat{a}+\hat{b}} J_{ab} = J_{ab}.   \label{Jab}
\eeq
This simplifies many relations in the later computations. 
Following Ref. \cite{JS} we express the elements 
$ e, \beta $ and $ \gamma $ in terms of the remaining
elements $ a, b, c, d, \alpha $. The commutation  relations  
satisfied by the elements $ a, b, c, d, \alpha $ 
and $ \delta $ are summarized as
\beq
  \begin{array}{lclcl}
    \vspace{1mm}
    [a, b] = h(1-a^2), & \ & [a,c] = hc^2, & \  & [a,d] = h(cd-ca), \\ 
    \vspace{1mm}
    {[a, \alpha]} = 0, & & {[a, \delta]} = hc\delta, & & {[b,c]} = h(ca+dc), \\
    \vspace{1mm}
    {[b, d]} = h(d^2-1), & & {[b, \alpha]} = h\alpha a, & & 
    {[b, \delta]} = h(d\delta + c \alpha), \\ \vspace{1mm}
    {[c, d]} = -hc^2, & & {[c, \alpha]} = -hc\delta, & & {[c, \delta]} = 0, \\ 
    \vspace{1mm}
    {[d, \alpha]} = h(\delta a - \delta d), & & {[d, \delta]} = h\delta c, & & 
    {\{ \alpha, \delta \}} = h(ac - \delta^2),\\
    \alpha^2 = {\displaystyle \frac{h}{2}(a^2-1)}, & & 
    \delta^2 = {\displaystyle \frac{h}{2}c^2}. & & 
  \end{array}
   \label{OSPcomm}
\eeq
The other entries of $T$ may be algebraically solved as follows:
\beq
  e = 1 + \alpha \delta - \frac{h}{2}ac, \qquad
  \beta = \alpha d - \delta b - h \delta d - \frac{h}{2} \gamma, \qquad 
  \gamma = \alpha c - \delta a -h \delta c. 
  \label{ebg}
\eeq
Relations analogous to the classical supergroup $OSp(2/1)$ exist for the
nonstandard deformation:
\beq
  ad - bc + \alpha \delta + \frac{h}{2}ac = 1, \qquad
  e^{-1} = (1- \alpha \delta + \frac{h}{2}ac) (1-\frac{h^2}{4}c^2)^{-1}, \qquad 
  \alpha \delta + \beta \gamma = \frac{h}{2}(ac-dc).
  \label{center}
\eeq
For completeness, we also give the commutation relations involving the
elements $ e, \beta $ and $ \gamma $:
\beq
  \begin{array}{lclcl}
     \vspace{1mm}
     [a, e] = h\gamma \delta, & & [b, e] = h(\beta \delta + \gamma \alpha), & & 
     [c, e] = 0, \\ \vspace{1mm}
     [d, e] = h \gamma \delta, & & [e, \alpha] = h(e\delta + \gamma a), & & 
     [e, \beta] = h(d \delta + \gamma e), \\ \vspace{1mm}
     [e, \gamma] = h c \delta, & & [e, \delta] = h c \gamma, & & 
     [a, \beta] = h (\gamma d - \gamma a), \\ \vspace{1mm}
     [b, \beta] = h \beta d, & & [c, \beta] = -h c \gamma, & & 
     [d, \beta] = 0, \\ \vspace{1mm}
     \{ \alpha, \beta \} = h(ea - ed), & & 
     \{ \beta, \gamma \} = -h(dc + \gamma^2), & & 
     \{ \beta, \delta \} = hce, \\ \vspace{1mm}
     [a, \gamma] = h\gamma c, & & [b, \gamma] = h (\beta c + \gamma a), & & 
     [c, \gamma] = 0, \\ \vspace{1mm}
     [d, \gamma] = h c\gamma, & & \{ \alpha, \gamma \} = -hce, & & 
     \{ \gamma, \delta \} = 0, \\
     \beta^2 = {\displaystyle \frac{h}{2}}(1-d^2), & & 
     \gamma^2 = -{\displaystyle \frac{h}{2}}c^2. & & 
  \end{array}
  \label{OSPcommadd}
\eeq
As a consequence of the grading the RTT relation involves extra sign 
factors in the tensor products of $T$ and the identity matrix:
\beq
  \sum_{x,y} (-1)^{\h{y}(\h{x}+\h{i})} R^{k \ell}_{xy}\, t^x_i\, t^y_j = 
  \sum_{x,y} (-1)^{\h{y}(\h{k}+\h{x})} t^{\ell}_y\, t^k_x\, R^{xy}_{ij}.
  \label{RTT}
\eeq

  The coalgebra mappings of the quantum supergroup $OSp_h(2/1)$ are, as 
usual, given by
\beq
  \Delta(T) = T \stackrel{\cdot}\otimes T, \qquad \epsilon(T) = 
\hbox{diag}(1,1,1).
  \label{coprocouni}
\eeq
The antipode is obtained from the coproduct:
\beq
  S(T) = \left(
    \begin{array}{ccc}
      d + \frac{h}{2}c & -\beta-\frac{h}{2}\gamma & 
      -b - \frac{h}{2}(a-d) + \frac{h^2}{4}c \\
      \delta & e & -\alpha+ \frac{h}{2}\delta \\
      -c & \gamma & a-\frac{h}{2}c
    \end{array}
  \right) 
  = 
  J^{-1} T^{st} J. \label{atpOSP}
\eeq
It is easy to see that $ T S(T) = S(T) T = \hbox{diag}(1,1,1).$ 

%
%
%
\setcounter{equation}{0}
\section{Differential calculus on quantum superspace}
\label{DiffCal}

 In this Section, a quantum superspace covariant under the action of 
$OSp_h(2/1)$ is introduced and a differential calculus in the sense of Wess 
and Zumino \cite{WZ} is constructed. The quantum superspace is a graded 
algebra, denoted by $\A$, generated by two odd ($\theta_1, \theta_2$) 
and one even ($x$) elements. The defining relations of the algebra 
$ \A $ read 
\bea
  & & [\theta_1, x] = -hx \theta_2, \qquad \{ \theta_1, \theta_2 \} = 0, 
      \qquad 
      [\theta_2, x] = 0,
  \nn \\
  & & \theta_1^2 = -\frac{h}{2}(x^2 - 2 \theta_1 \theta_2), \qquad 
      \theta_2^2 = 0.  \label{qplane}
\eea
It is straightforward to verify that the relations (\ref{qplane}) are 
preserved under the action of $OSp_h(2/1)$ from the left 
\beq
  \left(
    \begin{array}{c}
       \theta'_1 \\  x' \\  \theta'_2 
    \end{array}
  \right) 
  = 
  \left(
   \begin{array}{ccc}
      a & \alpha & b \\
      \gamma & e & \beta \\
      c & \delta & d
   \end{array}
  \right)   
  \left(
    \begin{array}{c}
       \theta_1 \\  x \\  \theta_2 
    \end{array}
  \right).
  \label{leftaction}
\eeq
The quantum superspace (\ref{qplane}) has an important difference 
from that associated to the Jordanian quantum supergroup $ GL_h(1/1) $ 
discussed in Refs. \cite{DP,C}. In the quantum superspace covariant 
under the action of $GL_h(1/1)$ , the deformation parameter $h$ is 
a Grassmann variable, whereas in (\ref{qplane}) the quantity $h$ commute 
with all elements of the quantum superspace. 

  A scalar element $\varphi$ exists in the quantum superspace $\A$:
\beq
  \varphi \equiv X^{st} J X = x^2 - 2\theta_1 \theta_2, \label{invzero}
\eeq
where
\beq
  X = \left(
    \begin{array}{c}
       \theta_1 \\  x \\  \theta_2 
    \end{array}
  \right),
  \qquad
  X^{st} = (-\theta_1, x, -\theta_2).
  \label{XJX}
\eeq
Then it is easy to show that $ \varphi $ is preserved by the left 
action of $OSp_h(2/1)$ Note that the parity of components of $ X $ is 
$ \h{X}^i = 1 + \h{i} $ (mod 2). The fourth relation in (\ref{qplane}) 
implies that $ \theta_1^2 $ is also a scalar in the algebra $\A. $ 
Employing a solution of the nongraded Yang-Baxter equation the defining 
relations (\ref{qplane})  of the algebra $ \A $ may be 
written in a compact form: 
\beq
  X^i X^j = \sum_{k,\ell} B^{ij}_{k\ell} X^{\ell} X^k, \qquad 
  B_{12} B_{13} B_{23} = B_{23} B_{13} B_{12},  \label{qplanR}
\eeq
where the matrix $B$ reads
\beq
   B(h) = \left(
     \begin{array}{ccc|ccc|ccc}
       -1 & \cdot & -h & \cdot & -h &  \cdot & h & \cdot & -h^2/2 \\
        \cdot & 1 &  \cdot & \cdot &  \cdot & -h & \cdot & \cdot & \cdot \\
        \cdot & \cdot & -1 & \cdot &  \cdot &  \cdot & \cdot & \cdot & -h \\ 
        \hline
        \cdot & \cdot & \cdot  & 1 &  \cdot &  \cdot & \cdot & h & \cdot \\
        \cdot & \cdot & \cdot  & \cdot &  1 &  \cdot & \cdot & \cdot & -h \\
        \cdot & \cdot & \cdot  & \cdot &  \cdot &  1 & \cdot & \cdot & \cdot \\
        \hline
        \cdot & \cdot & \cdot  & \cdot &  \cdot & \cdot  & -1 & \cdot & h \\
        \cdot & \cdot & \cdot  & \cdot &  \cdot & \cdot  & \cdot & 1 & \cdot \\
        \cdot & \cdot & \cdot  & \cdot &  \cdot & \cdot  & \cdot & \cdot & -1
     \end{array}
   \right)_.
   \label{Bmatrix}
\eeq
The matrix $B$ is related to the $R$ matrix of the deformed 
supergroup $OSp_h(2/1)$ as 
\beq
   R(h)^{k\ell}_{xy} = (-1)^{1+\h{k}+(1+\h{x})\h{y}}(B(h)^{-1})^{k\ell}_{xy}
   = (-1)^{1+\h{k}+\h{y}+\h{k}\h{l}} B(-h)^{k\ell}_{xy}, \label{RtoB}
\eeq
where the last equality follows from the relation 
$ (B(h)^{-1})^{k\ell}_{xy} = (-1)^{\h{k}\h{\ell}+\h{x}\h{y}} 
B(-h)^{k\ell}_{xy}. $ 

  The differential calculus on quantum space is an algebra generated by 
coordinates $X^i$, differentials $ \Xi^i \equiv dX^i$ and derivatives 
$ \partial_i = \frac{\partial}{\partial X^i}. $ 
The parity of differentials and derivatives are, in general, 
$ \h{\Xi}^i = 1 + \h{X}^i, \ \h{\partial}_i = \h{X}^i. $ 
The differential calculus on quantum superspace using solutions of the 
nongraded Yang-Baxter equation is developed in Ref. \cite{KU}. 
These authors require the exterior derivative $d$, which 
maps a $k$-form to a $(k+1)$-form, to maintain three properties: 
(i) nilpotency, (ii) graded Leibnitz rule
\beq
  d(f \wedge g) = (df) \wedge g + (-1)^{\hat{f}}f \wedge dg, \qquad
  f \in \Omega^p, \ g \in \Omega^q  
  \label{dLeib}
\eeq
and (iii) its action on a function $f(X^i)$ is given by 
$ df = {\displaystyle \sum_{i} \Xi^i \p{i} f}. $ 
Employing these properties, the following commutation relations among 
$ X^i,\,\Xi^i $ and $\p{i}$ may be determined: 
\bea
  & & \Xi^i \wedge \Xi^j = \sum_{k,\ell} (-1)^{\h{X}^i+\h{\Xi}^{\ell}} 
      B^{ij}_{k\ell} \, 
      \Xi^{\ell} \wedge \Xi^k,  \qquad
      X^i \Xi^j = \sum_{k,\ell} (-1)^{\h{X}^i} B^{ij}_{k\ell}\, 
      \Xi^{\ell}\, X^k,
      \nn \\
  & & \partial_j X^i = \delta_{ij} + \sum_{k,\ell} B^{i\ell}_{kj}\, 
      X^k \partial_{\ell},
      \hspace{2.5cm}
      \partial_j\, \Xi^i = \sum_{k,\ell} (-1)^{\h{X}^j}(B^{-1})^{ki}_{j\ell}\, 
      \Xi^{\ell}\, \partial_k,
      \nn \\
  & & \partial_i \partial_j = \sum_{k,\ell} B^{k\ell}_{ij}\, \partial_{\ell} 
      \partial_k.
      \label{DC} 
\eea
Our convention of the matrix $B$ differs from that in Ref. \cite{KU}. 
We use the Yang-Baxter equation of the form (\ref{qplanR}), whereas the 
Yang-Baxter equation in the braid group form 
$ F_{12} F_{23} F_{12} = F_{23} F_{12} F_{23} $ is used in Ref. \cite{KU}. 
They are related as $ F^{k\ell}_{ij} = B^{k\ell}_{ji}.$ 
The relations (\ref{DC}) are covariant under the  action of 
the super-Jordanian $OSp_h(2/1)$:
\beq
  X'^i = \sum_j t^i_j \, X^j, \qquad 
  \Xi'^i = \sum_j (-1)^{\h{i}+\h{j}}\, t^i_j \, \Xi^j, \qquad
  \partial'_i = \sum_j (-1)^{\h{i}+\h{j}} ((T^{st})^{-1})^i_j \, 
  \partial_j.
  \label{trans}
\eeq
To show the covariance, we need RTT type relations for $T$ and $T^{st}$ 
with the matrix $B$. They are obtained {\it via} (\ref{RTT}) and 
(\ref{qplanR}):
\bea
  & & \sum_{i,j} (-1)^{\h{b}+\h{j} + \h{i}\h{j} + \h{b}\h{i}}\, 
      t^a_i\, t^b_j\, B^{ij}_{cd}
      = \sum_{i,j} (-1)^{\h{c}+\h{i}+\h{c}\h{d}+\h{d}\h{i}}\, 
      B^{ab}_{ij}\, t^j_d\, t^i_c,
    \label{RTT1} \\
  & & \sum_{i,j} (-1)^{\h{i}+\h{a}+\h{i}\h{j}+\h{i}\h{b}}\, 
      \tau^a_i\, \tau^b_j \, B^{cd}_{ij} 
      = \sum_{i,j} (-1)^{\h{d}+\h{j}+\h{d}\h{i} + \h{c}\h{d}}\, 
      B^{ij}_{ab}\, 
      \tau^j_d\, \tau^i_c, 
    \label{RTT2} \\
  & & \sum_{i,j} (-1)^{\h{i}+\h{b}\h{i}+\h{c}\h{i}+\h{d}\h{j}}\, 
     \tau^a_i\, t^b_j\, B^{jd}_{ci} 
      = \sum_{i,j} (-1)^{\h{c}+\h{i}+\h{j}+\h{i}\h{c}}\, B^{bi}_{ja}\, 
      t^j_c\, \tau^i_d,
    \label{RTT3} \\
  & & \sum_{i,j} (-1)^{\h{c}+\h{i}+\h{j}+\h{c}\h{i}}\, 
     (B^{-1})^{ib}_{aj}\, t^j_c\, \tau^i_d 
     = \sum_{i,j} (-1)^{\h{i}+\h{b}\h{i}+\h{c}\h{i}+\h{j}\h{d}}\, 
     \tau^a_i\, t^b_j\, (B^{-1})^{dj}_{ic},
    \label{RTT4}
\eea
where $ \tau = (T^{st})^{-1}$. Introducing the notations 
$ \xi_1 = d\theta_1,\ \eta = dx, \ \xi_2 = d\theta_2, $ 
the explicit form of the $OSp_h(2/1)$ covariant differential calculus 
on the quantum superspace $ \A $ is summarized as follows:

\begin{itemize}
  \item coordinates
  \bea
  & & [\theta_1, x] = -hx \theta_2, \qquad \{ \theta_1, \theta_2 \} = 0, \qquad 
      [\theta_2, x] = 0,
  \nn \\
  & & \theta_1^2 = -\frac{h}{2}(x^2 - 2 \theta_1 \theta_2), \qquad 
      \theta_2^2 = 0.  \label{coordinates}
  \eea
  \item differentials 
  \bea
   & & \xi_1 \wedge \eta - \eta \wedge \xi_1 = h \eta \wedge \xi_2, \qquad
       \xi_1 \wedge \xi_2 - \xi_2 \wedge \xi_1 = h\xi_2 \wedge \xi_2,
   \nn \\
   & & \eta \wedge \xi_2 - \xi_2 \wedge \eta = 0, \qquad \qquad \ \;
       \eta \wedge \eta = -\hot \xi_2 \wedge \xi_2.  \label{differentials}
  \eea
  \item coordinates and differentials
  \bea
    & & [\theta_1, \xi_1] = h( \theta_1 \xi_2  + x \eta - \theta_2 \xi_1 - 
        \frac{h}{2} \theta_2 \xi_2), 
    \qquad
        \{ \theta_1, \eta \} = hx \xi_2, 
    \nn \\
    & & [\theta_1, \xi_2] = h\theta_2 \xi_2,
    \qquad
        [x, \xi_1] = -h \theta_2 \eta, 
    \qquad \quad
        [x, \eta] = -h \theta_2 \xi_2,
    \nn \\
    & & [x, \xi_2] = 0, \qquad \qquad \ \ 
        [\theta_2, \xi_1] = -h \theta_2 \xi_2,
     \label{cordiff} \\
    & & \{ \theta_2, \eta \} = 0, \qquad \qquad \ 
        [\theta_2, \xi_2] = 0. \nn
  \eea
  \item derivatives and coordinates
  \bea
    & & \p{1} \theta_1 = 1 - \theta_1 \p{1} + h \theta_2 \p{1}, \qquad
        \p{1} x = x \p{1}, \hspace{3cm} \p{1} \theta_2 = -\theta_2 \p{1},
    \nn \\
    & & \p{x} \theta_1 = \theta_1 \p{x} - h x\p{1}, \qquad\qquad
        \p{x} x = 1 + x\p{x} + h \theta_2 \p{1}, \qquad
        \p{x} \theta_2 = \theta_2 \p{x},
    \nn \\
    & & \p{2} \theta_1 = -\theta_1 \p{2} - h(\theta_1 \p{1} + x \p{x}
        + \theta_2 \p{2} + \frac{h}{2} \theta_2 \p{1}),
    \label{dercor} \\
    & & \p{2} x = x \p{2} - h\theta_2 \p{x}, \qquad \qquad
        \p{2} \theta_2 = 1 - \theta_2 \p{2} + h \theta_2 \p{1},	
    \nn
  \eea
  where $ \p{1} = {\displaystyle \frac{\partial}{\partial \theta_1}}, \ 
  \p{x} =  {\displaystyle \frac{\partial}{\partial x}}, \ 
  \p{2} =  {\displaystyle \frac{\partial}{\partial \theta_2}}.
  $
  \item derivatives and differentials
  \bea
    & & \p{1} \xi_1 = \xi_1 \p{1} - h \xi_2 \p{1}, \qquad\quad
        \p{1} \eta = -\eta \p{1}, \qquad\qquad\quad\ 
        \p{1} \xi_2 = \xi_2 \p{1},
    \nn \\
    & & \p{x} \xi_1 = \xi_1 \p{x} - h \eta \p{1}, \quad\qquad\;
        \p{x} \eta = \eta \p{x} + h \xi_2 \p{1}, \qquad
        \p{x} \xi_2 = \xi_2 \p{x},
    \nn \\
    & & \p{2} \xi_1 = \xi_1 \p{2} + h( \xi_1 \p{1} + \eta \p{x} 
        + \xi_2 \p{2} + \frac{h}{2} \xi_2 \p{1}),
    \label{derdiff} \\
    & & \p{2} \eta = -\eta \p{2} + h \xi_2 \p{x}, \qquad\quad
        \p{2} \xi_2 = \xi_2 \p{2} - h \xi_2 \p{1}.
    \nn
  \eea
  \item derivatives
  \bea
    & & \p{1}^2 = 0, \qquad\quad \p{1}\p{x} = \p{x}\p{1}, \qquad \quad
        \p{1} \p{2} = - \p{2}\p{1},
    \nn \\
    & & \p{x}\p{2} = \p{2} \p{x} - h \p{1} \p{x}, \qquad\qquad\quad\ \;
        \p{2}^2 = h(\p{1}\p{2} - \frac{1}{2}\p{x}^2).
    \label{derder}
  \eea
\end{itemize}

%
%
%
\setcounter{equation}{0}
\section{$OSp_h(2/1)$ symmetric torsionless connections}
\label{Connection}

  We have seen that a scalar $\varphi\,(\sim \theta_1^2) $ exists in the 
quantum superspace $ \A. $ This scalar is an $OSp_h(2/1)$ invariant 
zero-form. 
Invariant one and two-forms under the action of the deformed supergroup
$OSp_h(2/1)$ also exist in the differential calculus $\A$: 
\bea
  & & \varrho = \sum_{a,b} J_{ab}\, X^a \Xi^b = \theta_1 \xi_2 + x \eta 
      - \theta_2 \xi_1 - \hot \theta_2 \xi_2,
      \label{invone} \\
  & & \chi = \sum_{a,b} J_{ab}\, \Xi^a \wedge \Xi^b = 0. \label{invtwo}
\eea
It is evident that the invariant two-form $\chi$ is trivial.
It is straightforward to verify the invariance of $ \varrho $ and $\chi$ 
under the transformation (\ref{trans}). Note that the $ \varrho $ appears 
on the right hand side of the first relation in (\ref{cordiff}). It is 
easy to find the commutation relations between the invariant forms and 
the basis elements $(X^a, \Xi^a)$. For the zero-form $ \varphi $ these 
relations read 
\beq
    X^a \varphi = \varphi X^a, \qquad \Xi^a \varphi = \varphi \Xi^a.
    \label{XPhi}
\eeq 
The commutation properties of the invariant one-form $ \varrho $
are succinctly given by
 \beq
   X^a \varrho = (-1)^{\hat{X}^a} \varrho X^a, 
   \qquad
   \Xi^a \wedge \varrho = (-1)^{\hat{\Xi}^a} \varrho \wedge \Xi^a.
   \label{XXirho}
 \eeq
In a more expanded version the above relations read
 \bea
   & & [x, \varrho] = \{ \theta_i, \varrho \} = 0, \qquad i = 1, 2 
   \label{Xrho} \\
   & & \eta \wedge \varrho + \varrho \wedge \eta = 0, \qquad 
       \xi_i \wedge \varrho - \varrho \wedge \xi_i = 0. \label{Xirho}
 \eea
 It is also straightforward to verify the relation
 \beq
   \varrho \wedge \varrho = 0. \label{rho2}
 \eeq

   In order to determine the covariant derivative, it is necessary to 
find the action of the extended permutation $\sigma$ on 
$ \Omega^1 \otimes \Omega^1. $ This can be done by applying the 
covariant derivative $D$ on the second relation in (\ref{DC}). Using
the Leibnitz rules, we obtain
\[
   \Xi^i \otimes \Xi^j + (-1)^{\hat{X}^i} X^i D\Xi^j 
   =  \sum_{k,\ell} (-1)^{\hat{X}^i} B^{ij}_{k\ell}\; 
   \{ (-1)^{\hat{\Xi}^{\ell}} \sigma(\Xi^{\ell} \otimes \Xi^k) 
   + (D \Xi^{\ell}) X^k \}.
\]
This relation implies that the action of $ \sigma $ on $ \Xi^{\ell} 
\otimes \Xi^k $, and the commutation relations between $ X^i $ and 
$ D\Xi^j$ may be consistently described as 
\bea
  & & \Xi^i \otimes \Xi^j = \sum_{k,\ell} (-1)^{\hat{X}^i+\hat{\Xi}^{\ell}} 
      B^{ij}_{k\ell}\;\sigma(\Xi^{\ell} \otimes \Xi^k), \label{Detsigma}\\
  & & X^i D\Xi^j = \sum_{k,\ell} B^{ij}_{k\ell}\; (D\Xi^{\ell}) X^k. 
      \label{DetD}
\eea
Using the property (\ref{RtoB}) the exchange relation (\ref{Detsigma}) may 
be solved yielding the action of $\sigma$ on the tensored space of
one-forms as follows:  
\beq
    \sigma(\Xi^k \otimes \Xi^{\ell}) = 
    \sum_{i,j} (-1)^{\hat{i}\hat{j}} R^{\ell k}_{ij}\; \Xi^i \otimes \Xi^j
    \equiv \sum_{i,j} \check{R}^{k \ell}_{ij}\; \Xi^i \otimes \Xi^j,
    \label{sigma1}
\eeq
The matrix $ \check R $ has two important properties, namely, $\check R$ 
is idempotent and satisfies a non-graded Yang-Baxter equation
\beq
    \check R^2 = 1, \qquad \check R_{12} \check R_{23} \check R_{12} 
    = \check R_{23} \check R_{12} \check R_{23}. \label{Rcheck}
\eeq
As a consequence of the exchange of the superscripts in the definition 
(\ref{sigma1}) of $\check R$, it satisfies a different form of 
Yang-Baxter equation from the one obeyed by $R$. The operator 
$\sigma$, therefore, follows identical properties:
\beq
    \sigma^2 = 1, \qquad
    \sigma_{12} \sigma_{23} \sigma_{12} = \sigma_{23} \sigma_{12} \sigma_{23}.
    \label{sigmapro}
\eeq
The map $ \sigma $ may now be explicitly written as follows
\bea
    & & \sigma(\xi_1 \otimes \xi_1) = \xi_1 \otimes \xi_1 - h( \xi_1 
        \otimes \xi_2 + \eta \otimes \eta - \xi_2 \otimes \xi_1 
        - \hot \xi_2 \otimes \xi_2), \nn \\
    & & \sigma(\xi_1 \otimes \eta ) = \eta \otimes \xi_1 + h \xi_2 
        \otimes \eta, \nn \\
    & & \sigma(\xi_1 \otimes \xi_2) = \xi_2 \otimes \xi_1 + h \xi_2 
        \otimes \xi_2,\nn \\
    & & \sigma(\eta \otimes \xi_1) = \xi_1 \otimes \eta - h \eta 
        \otimes \xi_2, \nn \\
    & & \sigma(\eta \otimes \eta) = -\eta \otimes \eta - h\xi_2 
        \otimes \xi_2,
    \label{sigmaaction} \\
    & & \sigma(\eta \otimes \xi_2) = \xi_2 \otimes \eta,
    \nn \\
    & & \sigma(\xi_2 \otimes \xi_1) = \xi_1 \otimes \xi_2 - h\xi_2 
        \otimes \xi_2, \nn \\
    & & \sigma(\xi_2 \otimes \eta) = \eta \otimes \xi_2, \nn \\
    & & \sigma(\xi_2 \otimes \xi_2) = \xi_2 \otimes \xi_2. \nn
\eea
The map $ \sigma $ being $\A$-bilinear, the following relations hold:
\bea
   & & \sigma(\xi_1 \otimes \varrho) = \varrho \otimes \xi_1, \qquad  
       \sigma(\eta \otimes \varrho) = - \varrho \otimes \eta, \qquad 
       \sigma(\xi_2 \otimes \varrho) = \varrho \otimes \xi_2,
   \nn \\
   & & \sigma(\varrho \otimes \xi_1) = \xi_1 \otimes \varrho, \qquad 
       \sigma(\varrho \otimes \eta) = -\eta \otimes \varrho, \qquad 
       \sigma(\varrho \otimes \xi_2) = \xi_2 \otimes \varrho,
    \label{sigmarho} \\
   & & \sigma(\varrho \otimes \varrho) = -\varrho \otimes \varrho. 
   \nn
\eea
The explicit form of the map $\sigma$ being known, the relations 
(\ref{sigmapro}) and (\ref{pisigma}) may be verified by direct
computation. 

  To derive the action of the covariant derivative $ D$ on $ \Xi^a,$
we compare the relation (\ref{DetD}) with (\ref{qplanR}) and the
second relation in (\ref{DC}). The comparison suggests that 
$ D\Xi^a $ contains $ X^a $ and $ \Xi^a $ as factors. Another 
important observation is that $ D\Xi^a $ has the same transformation 
property as $ X^a$ under the action of $OSp_h(2/1)$ namely, 
\beq
      D \Xi'^a = \sum_{i} t^a_i \, D \Xi^i. \label{OsponDXi}
\eeq
Thus the most general form of $ D\Xi^a $ may be given by
\beq
     D \Xi^a = c_0 X^a \varpi + c_1 (-1)^{\hat{a}} \Xi^a \otimes \varrho 
       + c_2 \varrho \otimes \Xi^a,
     \label{DXi1}
\eeq
where $ c_i \ (i = 0, 1,2) $ are real parameters and 
$ \varpi \in \Omega^1 \otimes \Omega^1 $ satisfies
\beq
    \varpi' = \varpi, \qquad X^a \varpi = \varpi X^a.
    \label{varpipro} 
\eeq
It is not difficult to see that each term on the right hand side of 
(\ref{DXi1}) has the same transformation property as $ X^a $ under the 
action of the deformed supergroup $OSp_h(2/1)$ Furthermore, each term of 
(\ref{DXi1}) satisfies the same commutation relation as (\ref{DetD}). 
As we have seen in the beginning of this Section, the $OSp_h(2/1)$ 
invariant two-form $\chi$ is trivial so that the only possible choice for 
$ \varpi $ is given by
\beq
    \varpi = \varrho \otimes \varrho. \label{piform}
\eeq
In this way, we have seen that (\ref{DXi1}) and (\ref{piform})  
describe the most general linear connection. 

  Let us recall that our main interest is in torsionless connections, 
as the torsionfree condition is necessary for making the curvature left 
$ \A$-linear. We restrict the linear connection obtained above 
to be torsionfree: $\Theta \Xi^a = 0$. As the nilpotency of $d$ 
constrains $ d\Xi^a = 0 $, we obtain 
\bea
   & & \Theta \Xi^a = -\pi \circ D\Xi^a 
    = -c_0 \varrho \wedge \varrho - c_1 (-1)^{\hat{\Xi}^a} \Xi^a 
    \wedge \varrho - c_2 \varrho \wedge \Xi^a
    \nn \\
   & & \qquad
    = - (c_1 + c_2) \varrho \wedge \Xi^a = 0,
   \label{torsionlesscond}
\eea
where we have used the relations (\ref{XXirho}) and (\ref{rho2}).  The 
torsionfree condition thus requires $ c_2 = - c_1. $ Therefore, the 
general form of the $OSp_h(2/1)$ symmetric torsionless connections 
is given by the following two-parameter family: 
\beq
    D \Xi^a = c_0 X^a \varrho \otimes \varrho + c_1 ( 
     (-1)^{\hat{a}} \Xi^a \otimes \varrho - \varrho \otimes \Xi^a ).
    \label{DXi2}
\eeq
More explicitly these connections read
\bea
    & & D \xi_1 = c_0 \theta_1 \varrho \otimes \varrho 
        + c_1( \xi_1 \otimes \varrho - \varrho \otimes \xi_1), 
    \nn \\
    & & D \eta = c_0 x \varrho \otimes \varrho - c_1 ( \eta \otimes \varrho 
       + \varrho \otimes \eta),
    \label{DXi3} \\
    & & D \xi_2 = c_0 \theta_2 \varrho \otimes \varrho 
        + c_1 ( \xi_2 \otimes \varrho - \varrho \otimes \xi_2).
    \nn
\eea
For the torsionless connections (\ref{DXi2}), it is easy to see
\beq
  D \varrho = \sum_{a,b} J_{ab}\, \Xi^a \otimes \Xi^b + (c_0 \varphi - 2c_1) 
  \varrho \otimes \varrho.
  \label{Drho}
\eeq
Applying (\ref{pidef}, \ref{invtwo}, \ref{rho2}) it immediately follows 
that
\beq
    \pi(D \varrho) = 0.  \label{piDrho}
\eeq

%
%
%
\setcounter{equation}{0}
\section{Curvature and metric}
\label{Curv-Met}

  A two-parameter family of $OSp_h(2/1)$ symmetric torsionfree connections 
was obtained in the previous section. Since the generalized permutation 
operator $ \sigma $ satisfies the relation (\ref{pisigma}), the curvature 
computed from the connections are left $\A$-linear. Recall that curvatures 
are, in general, not right $\A$-linear. In the present case, however, the 
curvature is also right $\A$-linear. We exhibit this by  explicit 
computation. We also discuss the metric on the quantum superspace $ \A$. 
It, however, turns out that the connections are not compatible with the 
metric. 

  To obtain the curvature, we apply $ \pi_{12} D $ on (\ref{DXi2}). 
Each term is computed separately and listed below:
\bea
  & & \pi_{12}D (X^a \varrho \otimes \varrho)
      = \Xi^a \wedge \varrho \otimes \varrho 
      - \sum_{b,c} (-1)^{\hat{X}^a} 
      J_{bc}\, X^a \varrho  \wedge \Xi^b \otimes \Xi^c,
  \nn \\
  & & \pi_{12}D ((-1)^{\hat{a}} \Xi^a \otimes \varrho)  
      = \sum_{b,c} J_{bc}\, \Xi^a \wedge \Xi^b \otimes \Xi^c
      + (c_0 \varphi - 2c_1) \Xi^a \wedge \varrho \otimes \varrho,
  \nn \\
  & & \pi_{12}D (\varrho \otimes \Xi^a) 
      = -c_1 \Xi^a \wedge \varrho \otimes \varrho.
  \nn
\eea
Combining the above results, the curvature is obtained as follows:
\beq
  \pi_{12} D^2 \Xi^a = (c_0 - c_1^2 + c_0 c_1 \varphi) 
     \Xi^a \wedge \varrho \otimes \varrho
     + (c_0 (-1)^{\hat{a}} X^a \varrho + c_1 \Xi^a) \wedge \Lambda,
  \label{Curv1}
\eeq
where
\beq
  \Lambda = \sum_{a,b} J_{ab}\, \Xi^a \otimes \Xi^b 
  = \xi_1 \otimes \xi_2 + \eta \otimes \eta - \xi_2 \otimes \xi_1
    - \hot \xi_2 \otimes \xi_2.
  \label{Curv2}
\eeq
Note that $ \pi(\Lambda) = \chi = 0. $  
Expanding the first term in the right hand side of (\ref{Curv1}) as
\[
  \Xi^a \wedge \varrho \otimes \varrho
  =
  \sum_{b,c} (-1)^{\hat{X}^b} J_{bc}\, \Xi^a X^b \wedge \varrho \otimes \Xi^c,
\]
we express the curvature in terms of a two-form $ \omega $
\beq
   \pi_{12}D^2 \Xi^a = \sum_b \omega^a_b \otimes \Xi^b,   \label{CurvTwo1}
\eeq
where
\beq
   \omega^a_b = \sum_k J_{kb} \left\{
   (-1)^{\hat{k}}  \{
     c_0 (-1)^{\hat{a}} X^a \Xi^k - (c_0 -c_1^2 + c_0 c_1 \varphi) \Xi^a X^k 
    \} \wedge \varrho 
    + c_1  \Xi^a \wedge \Xi^k \right\}.
   \label{CurvTwo2}
\eeq

  We now prove that the curvature obtained above is right $\A$-linear. 
To this end, we note that the following relation may be established 
by direct computation:
\beq
   [X^a,\, \Lambda ] = 0. \label{XaLambda}
\eeq
Employing the second relation in (\ref{DC}), in conjunction with
the left $\A$-linearity of the curvature, we obtain 
\beq
   \pi_{12}D^2(\Xi^b X^a) = 
   \sum_{i,j} (-1)^{\hat{X}^i}(B^{-1})^{ab}_{ij} X^i \pi_{12}D^2 \Xi^j.
   \label{Curv3}
\eeq
Substituting (\ref{Curv1}) into (\ref{Curv3}), and then 
transferring $ X^i $ to the right {\it via} equations 
(\ref{DC}, \ref{XPhi}, \ref{XXirho}, \ref{XaLambda}),  
we demonstrate the intended result
 \beq
   \pi_{12}D^2(\Xi^b X^a) = (\pi_{12}D^2 \Xi^b) X^a, \label{RightLin}
 \eeq
establishing the right $\A$-linearity of the curvature. In the above 
computation we have used the fact that the matrices $ B^{ab}_{ij}$ 
and $ (B^{-1})^{ab}{ij}$ maintain the following relationship
regarding the parity of their indices: 
$ \hat a + \hat b = \hat i + \hat j$. 

  Let us now turn to the metric, which is considered as a bilinear map 
$ g : \Omega^1 \otimes_{\A} \Omega^1 \rightarrow \A. $ 
To completely determine the metric we need to know the action of the 
map $g$ on the basis elements of $\Omega^1 \otimes \Omega^1$. 
Setting $g^{ab} = g(\Xi^a \otimes \Xi^b)$, we require that $ g^{ab} $ 
to be invariant under the action of $OSp_h(2/1)$:
\bea
   g'{}^{ab} &\equiv& g(\Xi'{}^a \otimes \Xi'{}^b) 
   = \sum_{k,\ell} g((-1)^{\hat{a}+\hat{k}}\, t^a_k \, \Xi^k \otimes 
       (-1)^{\hat{b}+\hat{\ell}}\, t^b_{\ell} \, \Xi^{\ell})
   \nn \\
   &=&  \sum_{k,\ell} (-1)^{\hat{a}+\hat{b}+\hat{k}+\hat{\ell}}\, 
   t^a_k \, g^{k\ell} (t^{st})_{\ell}^b.
   \nn
\eea
The above result, in conjunction with the identity (\ref{Osymh}),  
immediately yields $ g'{}^{ab} = g^{ab} $, provided we choose  
$ g^{k\ell} = (-1)^{\hat{k}+\hat{\ell}} (J^{-1})_{k\ell} 
= (J^{-1})_{k\ell}. $ 
The $OSp_h(2/1)$ invariant metric, therefore, is given by
\beq
   g^{ab} = g(\Xi^a \otimes \Xi^b) = (J^{-1})_{ab} 
   = \left(
     \begin{array}{ccc}
       -\hot & 0 & -1 \\
        0    & 1 &  0 \\
        1    & 0 &  0
     \end{array}
     \right).
   \label{metric}
\eeq
Denoting the components of $ g^{-1} $ by $ g_{ab}, $ we note that 
the invariant one-form $ \varrho $ may be written in terms of the metric
\[
 \varrho = \sum_{a,b} g_{ab} X^a \Xi^b.
\]
The structure  of the metric (\ref{metric}) implies
\beq
d \circ g(\Xi^a \otimes \Xi^b) = 0.
\label{dgzero}
\eeq
The compatibility condition (\ref{Dtog}) now reads
\beq
   (1 \otimes g) \circ D(\Xi^a \otimes \Xi^b) = 0.
\label{dDontwo}
\eeq
 To compute the left hand side in (\ref{dDontwo}), we start by ordering 
the one-forms in the expression of $\varrho$ to the left: 
\beq
 \varrho = \xi_2 \theta_1 + \eta x - \xi_1 \theta_2 + \hot \xi_2 \theta_2 
  = \sum_{a,b} (-1)^{\hat{X}^b} J_{ab} \, \Xi^a X^b.
  \label{varrho2}
\eeq
We now readily obtain
\beq
   g(\varrho \otimes \Xi^a) = X^a, \qquad 
   g(\Xi^a \otimes \varrho) = (-1)^{\hat{X}^a}X^a.
   \label{grhoX}
\eeq
Following (\ref{Donn}) the action of the covariant derivative on $ 
\Omega^1 \otimes \Omega^1 $ is given as
\beq
   D(\Xi^a \otimes \Xi^b) = D \Xi^a \otimes \Xi^b + 
   (-1)^{\hat{\Xi}^a} \sigma_{12} (\Xi^a \otimes D \Xi^b).
   \label{DXiXi}
\eeq
Substituting (\ref{DXi2}) into (\ref{DXiXi}), we observe that, as a
consequence of the bilinearity of $g$, we may treat the first 
(proportional to $c_0$) and the second (proportional to $c_1$) 
terms in the right hand side of (\ref{DXi2}) separately. 
For the choice $ c_1 = 0$, we then obtain 
\[
    D(\Xi^a \otimes \Xi^b) = c_0(
    X^a \varrho \otimes \varrho \otimes \Xi^b + \varrho \otimes \Xi^a 
     \otimes \varrho X^b ),
\]
which, in turn, yields
\beq
   (1 \otimes g) \circ D(\Xi^a \otimes \Xi^b) = (-1)^{\hat{X}^a}2 c_0 
   \varrho X^a X^b.
   \label{1gDmu}
\eeq
For the alternate choice $ c_0 = 0, $ it follows that
\[
  D(\Xi^a \otimes \Xi^b) = c_1 \{
   (-1)^{\hat{a}} \Xi^a \otimes \varrho \otimes \Xi^b 
   - 2 \varrho \otimes \Xi^a \otimes \Xi^b + (-1)^{\hat{a}+\hat{b}}
   \sigma_{12} (\Xi^a \otimes \Xi^b \otimes \varrho ) \}.
\]
The right hand side in (\ref{dDontwo}) now reads
\bea
   & & (1 \otimes g) \circ D(\Xi^a \otimes \Xi^b) 
   \nn \\ 
   & & \qquad
   = c_1 \{ (-1)^{\hat{a}} \Xi^a X ^b - 2 g^{ab} \varrho 
   +(-1)^{\hat{a}+\hat{b}} (1 \otimes g) \circ \sigma_{12}
   (\Xi^a \otimes \Xi^b \otimes \varrho) \}.
   \label{1gDrho}
\eea
The last term is computed by using (\ref{sigmaaction}) and (\ref{grhoX}). 
The result is listed below:
\bea
  & & (1 \otimes g) \circ D(\xi_1 \otimes \xi_1) = 0, \nn \\
  & & (1 \otimes g) \circ D(\xi_1 \otimes \eta) 
      = c_1 (\xi_1 x + \eta \theta_1 - h \xi_2 x),
  \nn \\
  & & (1 \otimes g) \circ D(\xi_1 \otimes \xi_2) 
  =  c_1 (\eta x - \hot \xi_2 \theta_2 + \varrho),
  \nn \\
  & & (1 \otimes g) \circ D(\eta \otimes \xi_1) 
  =  - c_1 (\xi_1 x + \eta \theta_1 + h \eta \theta_2),
  \label{1gDrho2} \\
  & & (1 \otimes g) \circ D(\eta \otimes \eta) 
  = -  c_1 (2\eta x - h\xi_2 \theta_2 + 2 \varrho),
  \nn \\
  & & (1 \otimes g) \circ D(\eta \otimes \xi_2) 
  =  - c_1 (\eta \theta_2 + \xi_2 x),
  \nn \\
  & & (1 \otimes g) \circ D(\xi_2 \otimes \xi_1) 
  =  - c_1 (\eta x - \hot \xi_2 \theta_2  + \varrho),
  \nn \\
  & & (1 \otimes g) \circ D(\xi_2 \otimes \eta) 
  =  c_1 (\eta \theta_2 + \xi_2 x),
  \nn \\
  & & (1 \otimes g) \circ D(\xi_2 \otimes \xi_2) = 0. 
  \nn
\eea
Together with (\ref{1gDmu}), it has been shown that
$ (1 \otimes g) \circ D \neq 0 $, except for the trivial choice 
$ c_0 = c_1 = 0. $ Thus the covariant derivative $D$ is not 
compatible with the metric.

%
%
%
\setcounter{equation}{0}
\section{Concluding remarks}
\label{CR}

 In the present work we have studied noncommutative spaces, linear
connections, curvatures and metrics associated with the quantized 
supergroups. Our approach is a naive extension of the differential 
geometry developed in \cite{DVM,Mou,DVMMMcurv}. We have demonstrated 
that the ideas of these authors may be appropriately adapted to study  
the geometric objects related to the quantum supergroups. Specifically, 
we applied the extended differential geometry to the quantum superspace 
covariant under the quantum supergroup $OSp_h(2/1)$. We have seen that 
our particular example has a two-parameter family of $OSp_h(2/1)$ symmetric 
torsionfree connections. It turned out that the curvature of the 
connection was bilinear. The connection was, however,  not compatible 
with the metric. These properties are specific to our example. 
There could be other quantum superspace endowed with linear 
connections compatible with metric. 

  It may be of interest to recall the results related to the 
quantum spaces covariant under quantized $SL(2)$ groups; and compare 
them with the present results. It is well-known that $SL(2)$ admits two 
inequivalent deformations: the standard $q$-deformation and the 
Jordanian $h$-deformation. The quantum space for $q$-deformed $ SL(2)$ 
has a one-parameter family of torsionless linear connections and it has 
been shown that there can be no compatible metric \cite{DVMMM}, whereas 
the quantum space of $h$-deformed $SL(2)$ is more classical. It has a 
two-parameter family of torsionfree linear connections. A one-parameter 
subfamily of these connections is known to be compatible with 
a metric \cite{CMP}. On the other hand, the Lie superalgebra $osp(2/1)$ 
admits three inequivalent deformations \cite{JS2}. We are thus able to 
consider three deformations of the supergroup $OSp(2/1)$: $q$-deformation 
\cite{KR}, $h$-deformation \cite{CP} and super-Jordanian deformation. 
The $q$ and $h$-deformations have the $ SL(2)$ counterparts, while 
super-Jordanian does not. The super-Jordanian deformation 
can be regarded as an algebra intermediate between $q$ and 
$h$-deformations. We have seen that the quantum space for super-Jordanian 
$OSp_h(2/1)$ is less classical since the connections are not metric. This 
leads us to the anticipate that the quantum space for $h$-deformed 
$OSp(2/1)$ has connections which are metric, while the connections 
on the quantum spaces related to the standard $q$-deformed supergroup 
$OSp_q(2/1)$  are not metric. This will be presented in a future work. 

%
%
\noindent{\large {\bf Acknowledgments}}\\
The work of N.A. was partially supported by the grants-in-aid from the 
MEXT, Japan. The other author (R.C.) was partially supported by the grant
DAE/2001/37/12/BRNS, Government of India.

%
%

\end{document}